\documentclass[11pt]{article}

\usepackage{mathrsfs} 

\usepackage[all]{xy}

\usepackage{color}

\definecolor{darkgreen}{rgb}{0.0, 0.7, 0.0}

\definecolor{cyan}{cmyk}{1,0,0,0}

\newtheorem{tm}{Theorem}[subsection]
\newtheorem{lm}[tm]{Lemma}
\newtheorem{pr}[tm]{Proposition}
\newtheorem{rmk}[tm]{Remark}

\newtheorem{ex}[tm]{Example}
\newtheorem{fact}[tm]{Fact}
\newtheorem{??}[tm]{Question}
\newtheorem{defi}[tm]{Definition}

\font\tenmsb=msbm10
\font\sevenmsb=msbm7
\font\fivemsb=msbm5
    
\newfam\msbfam
\textfont\msbfam=\tenmsb
\scriptfont\msbfam=\sevenmsb
\scriptscriptfont\msbfam=\fivemsb
\def\Bbb#1{{\fam\msbfam #1}}

\font\teneufm=eufm10
\font\seveneufm=eufm7
\font\fiveeufm=eufm5
\newfam\eufmfam
\textfont\eufmfam=\teneufm
\scriptfont\eufmfam=\seveneufm
\scriptscriptfont\eufmfam=\fiveeufm
\def\frak#1{{\fam\eufmfam\relax#1}}

\def\lorw{\longrightarrow}
\newcommand\n{\noindent}

\newcommand\ci{\cite}
\newcommand\s{\sigma}
\newcommand\rat{{\Bbb Q}}
\newcommand\comp{{\Bbb C}}

\newcommand\zed{{\Bbb Z}}

\newcommand\blacksquare{{\hspace*{\fill} $\fbox{}$}}

\newcommand\e{\epsilon}

\newcommand{\im}{ \hbox{\rm Im} }

\newcommand{\ke}{ \hbox{\rm Ker} }

\newcommand{\ben}{\begin{enumerate}}
\newcommand{\een}{\end{enumerate}}

\newcommand{\bit}{\begin{itemize}}
\newcommand{\eit}{\end{itemize}}

\newcommand{\beq}{\begin{equation}}
\newcommand{\eeq}{\end{equation}}

\newcommand{\la}{\label}

\newcommand{\m}[1]{\mathcal{#1}}
\newcommand{\ms}[1]{\mathscr{#1}}
\newcommand{\bb}[1]{\Bbb{#1}}

\title{The projectors of the decomposition theorem are motivic}
\author{
Mark Andrea A.  de Cataldo\thanks{
Partially supported by N.S.A. and N.S.F.}\, 
and Luca Migliorini\thanks{Partially supported by PRIN project ``Spazi di moduli e teoria di Lie''}}
\date{}

\begin{document}\maketitle

\begin{abstract}
We prove that the projectors arising from the decomposition theorem applied to a projective
map of quasi projective  varieties are absolute Hodge,   Andr\'e motivated,
Tate and   Ogus classes.  As a by-product,
we introduce,  in characteristic zero, the  notions of algebraic de Rham intersection cohomology groups
of a quasi projective  variety  and of  intersection cohomology motive
of a projective variety.
  \end{abstract}

\tableofcontents

\section{Introduction and preliminary material}\la{intro}
P. Deligne introduced the notion of absolute Hodge classes in [10], where he proved
that Hodge classes on complex Abelian varieties are absolute Hodge. This is a powerful
statement: on a complex Abelian variety, the notion of a Hodge class is purely algebraic.
Algebraic cycle classes are absolute Hodge,
and absolute Hodge classes are Hodge classes. A positive answer to the Hodge conjecture
would imply that these implications can be reversed. There is also a notion of absolute
Hodge map.

Let  $K$ be an algebraically closed  field of characteristic zero, let   $f: X \lorw Y$ be a projective map of  quasi projective
$K$-varieties  and   let 
$E$ be an $f$-ample line bundle on $X.$
In this paper, first we show how the decomposition  and relative hard Lefschetz theorems \ci{bbd,bams}
give rise to self-maps of the  intersection cohomology groups of $X,$ the {\em projectors of the
decomposition theorem}, and then we prove Theorem \ref{maintmz}: these projectors  are absolute Hodge maps. In particular, if $X$ is nonsingular projective,  then
these projectors are absolute Hodge classes on $X\times X.$
\S\ref{nomen} is devoted to prove some variants of Theorem
\ref{maintmz}: the projectors are absolute Hodge over any field of characteristic zero,
they are motivated in the sense of Y. Andr\'e, they are Tate classes and, finally, 
they are absolutely Hodge and Tate in the sense of A. Ogus.

There is no loss of generality in assuming that ${\rm tr.deg.}_\rat K < \infty.$
The notion of absolute Hodge involves the interplay of three
cohomology theories:
de Rham, \'etale $\rat_\ell$-adic and, after base change via an embedding of the field into $\comp,$
Betti.  

We take as starting point the decomposition and relative Hard Lefschetz theorems
in the Betti theory and, from that point on,  {\em we work exclusively in cohomology}, i.e. we make no further use
of  derived categories.  There are three reasons for this. The first is that we are not aware of
the existence in the literature of the cup-product operation at the level of
the derived category of D-modules (de Rham side). The second is that in the contexts
of motivated cycles and of crystalline cohomology, such derived techniques
do not seem to be available at the present time. The third is that we found working within the context
of cohomology and  of its fundamental functoriality properties aesthetically pleasing.

The key point in all the results of this paper, is to construct the projectors of the decomposition theorem  in a {\em uniform way in all  cohomology theories}, so that they turn out to be, more or less automatically, compatible
with each other via the comparison isomorphisms. 
We employ four main $K$-rational constructions from \S\ref{tgc}:  a geometric  construction  of the perverse filtration;  stratifications of maps; 
linear algebra description of  the decomposition by supports; splittings in abelian categories.
First, we work   with  $X$  nonsingular, then  we refine
 our analysis to the singular case.
Along the way, we offer in \S\ref{idrc} a definition of intersection de Rham cohomology
as a certain  subquotient of the de Rham cohomology of a resolution, and we  
point out several of its properties.

\smallskip
{\bf Acknowledgments.} We thank F. Charles and to B. Bhatt for useful 
conversations.

\subsection{Notation, conventions and some preliminary facts}\la{compbe}
 {\bf The set-up.}
Unless otherwise stated, we work over an algebraically closed
field $K$ of finite transcendence degree over $\rat$. 
 A $K$-scheme is a  separated scheme of finite type
  over $K;$ a $K$-variety is  
  an integral $K$-scheme. 
 The restriction on the transcendence degree of $K$ is for convenience of exposition only;
  see \S\ref{ahcz}.
The embeddings of $K$ into $\comp$ are denoted by $\s$. If $X$ is a $K$-scheme, we denote by
$\s X$ the pull-back to $\comp$; similarly for $K$-maps, etc.
 
  We employ the following cohomology theories (see \ci{dmos}): $H_{dR} (-/K)$ (de Rham),
  $H_{et}(-, \rat_\ell)$ ($\rat_\ell$-adic); if $K=\comp$, then we have $H_B(-,-)$ (the usual Betti 
  cohomology theory with $\rat, \comp, \rat_\ell$-coefficients. The compactly supported counterparts
  are denoted $H_{!,dR}(-/K)$, etc.
   We fix the following data: 
  \beq\la{data}
  \xymatrix{h := f\circ g : W \ar[r]^{\hskip 1.0cm g} & X \ar[r]^f &Y, &\eta := c_1(E),}
  \eeq
  where $f,g$ are projective maps of quasi projective $K$-schemes, $E$ is an $f$-ample
  line bundle on $X$ and $\eta:= c_1(E) \in H^2(X)(1),$ in any of the cohomology theories
  that we shall employ.
Our interest actually   lies in the datum of $(f,\eta).$ 
  The map $g$  is usually going to be a  resolution of the  singularities of $X,$
  and it is introduced as means  to reduce proofs of statements about $(f,\eta)$ to the case when 
  $X$ is nonsingular.

  {\bf Conventions on cohomological and filtration degrees}.
Since the bookkeeping of cohomological and filtration degrees  is inessential and a distraction, we omit it 
for the most part. Here is the  list  of the conventions that, unless mentioned otherwise, we employ implicitly: 
the cohomology  groups  $H^k(X)$ 
  live in the range $[0, 2\dim{X}],$ the Hodge filtrations
${\mathscr F}$
(Betti and de Rham)
in the range $[0, \dim{X}],$  and the weight  filtrations $\mathscr{W}$ (all three theories) in
the interval $[0, 2 \dim{X}];$
 if $X$ is integral and nonsingular,
then the  Betti intersection complex is $IC_X= \rat_X [\dim{X}];$
if $X$ is irreducible, then the cohomology sheaves ${\cal H}^k (IC_X)$
live in the range $[-\dim{X}, -1]$ and the    intersection cohomology
groups   $I\!H^k(X) := {\Bbb H}^{k-\dim X}(IC_X) $   live in the
same range as cohomology, i.e.  $[0, 2\dim{X}];$ 
for general $X,$ we write $X= \cup X_i$ (union of irreducible components),
form the finite map $\xymatrix{\nu: \coprod_i X_i  \ar[r] & X}$ and 
set $IC_X:= \nu_* ( \oplus_i IC_{X_i})$;  we place
the groups $I\!H^k(X_i)$ in the range $[0, 2 \dim{X_i}]$ 
(\ci{decpf2}, \S4.6);  the perverse filtration ${\mathscr P}_f$ on $I\!H(X)$ lives in the range
$[-r(f), r(f)],$ where $r(f)$ is a convenient integer, called the defect of semismallness of $f$
in \ci{htam, bams}.  We employ the same conventions for $H_!(X)$ and $I\!H_! (X)$. 

{\bf Decomposition, relative hard Lefschetz and semisimplicity theorems.}
We refer to  the survey \ci{bams} for the language and facts surrounding these
theorems, which we use freely in what follows, especially in \S\ref{dtbet}.

   {\bf Supports.} 
   A semisimple perverse sheaf $P$ (\'etale or Betti) on a $K$-scheme $V$ splits canonically
 $P\cong \oplus_{{\frak V} \in {\cal S}(P)} P_{\frak V}$ 
   as the  finite direct sum, over a finite set  $\m{S}(P)$ of distinct integral subvarieties ${\frak V} \subseteq V,$
of the  intersection complexes $P_{\frak V}$  of
the varieties  ${\frak V}$  with suitable twisted semisimple  coefficients.  We call $\m{S}(P)$ the {\em set of supports of $P.$} 
Let $\xymatrix{\phi: U \ar[r]& V}$ be a projective map of quasi-projective $K$-schemes, and let $Q$ be a semisimple
  complex 
  of geometric origin (\ci{bbd}, 6.2.4), e.g. $Q =IC_U.$ By the decomposition theorem, we have   a  non-canonical finite  direct sum decomposition:
$ \phi_* Q \cong \bigoplus_a Q_{\phi,a} [-a] $,   where   
$Q_{\phi,a}: = {^{\frak p}\!R}^a \phi_* Q$ 
 (perverse direct image sheaf).
Each $Q_{\phi,a}$ is a semisimple  perverse sheaf and we  obtain
 $\m{S} (\phi, Q, a) := \m{S} (Q_{\phi, a}),$
 the {\em set of supports of $\phi_* Q$ in perversity $a.$} The cohomology groups $H(U, Q)= H(V, \phi_*Q)$ are filtered by the perverse filtration
$\ms{P}_\phi$ and, given that $Q$ will be fixed by the context, we set
$Gr_{\phi,a}:= Gr^{\ms{P}_\phi}_a H(U,Q)$  and we denote the resulting
{\em decomposition by supports} by 
  $Gr_{\phi, a} = \oplus_{\frak{V} \in \m{S} (\phi, Q, a) } Gr_{\phi,a, \frak{V}}.$

{\bf Comparison Betti-de Rham and  Betti-\'etale.} 
 The algebraic de Rham cohomology groups
  $H_{dR}(X/K)$ are defined via  the same simplicial methods of \ci{ho3}.
They carry  the two natural filtrations
  $\mathscr{F}$ (Hodge) and $\mathscr{W}$ (weight). The usual maps between them, e.g. pull-backs via algebraic maps,
  are strict for these filtrations. 
  Recall that a filtered map is strict if taking graded objects is  exact. For any embedding $\s$
  of $K$ into $\comp$, 
there are the natural comparison isomorphisms
  $H_B(\s X, \comp) \cong H_{dR} (\s X/\comp)$ which are bi-filtered strict
  for $\mathscr{F}$ and $\mathscr{W}.$ Similarly, we have the natural comparison 
  isomorphisms  $H_B(\s X, \rat_\ell) \cong H_{et} (\s X, \rat_\ell)$    which is filtered
  strict for the weight filtrations on both sides. All of the above applies to $H_{!}$.

\subsection{Various notions of motivic endomorphism in cohomology}\la{vnmec}
\subsubsection{Diagram $(*)$ and the notion of absolute Hodge endomorphism}\la{gut}
  Let $(K,X,\s)$ be as in (\ref{data}).
  We have the following diagram expressing the relations among the various 
  cohomology groups $H(X),$
  Betti, algebraic de Rham, and $\rat_\ell$-adic, endowed with their
  Hodge and weight filtrations $\ms{F,W}$:
    \beq\la{pz1}
  \xymatrix{
   &  H_B (\s X, \comp) \ar[r]^{\hskip -0.2cm \cong}_{\hskip -0.2cm
   c_{ \mathscr{F},\mathscr{W}} } & H_{dR} (\s X/\comp) & \ar[l]_{-\otimes 1}^{\hskip 0.2cm 
   \s^*_{ \mathscr{F},\mathscr{W}}} H_{dR} ( X/K)
  \\
  {\hskip -2.5cm (*)} \qquad \quad H_B (\s X, \rat) \ar[ru]_{UCT_\mathscr{W}}^{-\otimes 1} \ar[rd]_{UCT_\mathscr{W}}^{-\otimes 1}
  \\
  & H_B (\s X, \rat_\ell) \ar[r]^{\hskip -0.2cm \cong}_{\hskip -0.3cm 
  c_{ \mathscr{W}} } & H_{et} (\s X, \rat_\ell) & \ar[l]_{\hskip 0.2cm \cong}^{\hskip 0.2cm 
  \s^*_{ \mathscr{W}} } H_{et} (X, \rat_\ell),
  }
  \eeq
  where: the arrows  UCT come from the universal coefficient theorem; the arrows $c$ are the comparison isomorphism for Betti-de Rham and Betti-\'etale cohomology; the arrows $\s^*$ come from
  base change via $\s;$ the arrows  are  strictly compatible for the indicated
  filtrations. 
  
   There is an analogous diagram for $H_!$  and  all the considerations
  and definitions 
  that follow apply to them. We leave them implicit.

\begin{defi}\la{defez}
  {\rm ({\bf Diagram $(*)$})
  We call  (\ref{pz1})  {\em diagram $(*)$ for $H(X).$} Such a diagram is made of vertices and arrows.
   Let $V$ be a direct sum of  tensor products of Tate-twisted subquotients
  of $H(X),$ obtained via some construction carried out at the various vertices. We say that we have {\em diagram $(*)$ for $V$} if, for every $\s,$  we can replace $H(X)$ with $V$ at all the vertices in (\ref{pz1}) and retain all the properties  listed above
 of the resulting    arrows induced by the construction. Similarly, given a map $\xymatrix{u: V \ar[r] & V'}$ between two such objects, obtained via some construction carried out at each vertex, we say that we have {\em diagram $(*)$ for $u$}
if we have diagrams $(*)$ for $V$ and $V'$ and the arrows $u$ yield a 
map of the  diagrams $(*)$ that makes all squares commutative and such all the arrows between
corresponding vertices of the two diagrams $(*)$ are filtered strict for the indicated filtrations.
}
\end{defi}

\begin{ex}\la{esempio}{\rm 
We have diagrams 
 $(*)$ for $V:= {\rm End} (H(X)), H(X)(m).$ 
 Proposition \ref{pfg} shows that we have  diagrams $(*)$ for the subspaces of the perverse filtration
 $\mathscr{P}_f$ 
 and hence for its graded pieces $Gr_f:= Gr^{\mathscr{P}_f}.$ Propositions \ref{hji}, \ref{rrtt}  show that we have diagrams
 $(*)$ for 
 the 
 direct sum decompositions  (\ref{rbb5}) and for the splitting (\ref{io44}).
}
\end{ex}

One may think of what above as saying that  a given construction $\xymatrix{u: V \ar[r] & V'},$
when seen at the de Rham vertex corresponding to $\s X,$
 is
$\rat$-rational  {\em as well as}  $K$-rational.

\bigskip
  Define:
  \beq\la{bny}
  \underline{H}(X)(m):= H_{dR}(X/K)(m) \times \prod_\ell H_{et}(X,\rat_\ell)(m).
  \eeq
  We define  $\underline{H}(\s X)(m)$  in a similar way.
  There are  natural maps ($\delta_\rat$: product of  comparison maps (\ref{pz1});
  $\s^*$: base change via $\s$):
  \beq\la{3ef}
   \xymatrix{H_B (\s X, \rat) (m)  \ar[r]^{\hskip 0.3cm \delta_\rat} & \underline{H}(\s X)(m)
   & \ar[l]_{\hskip 0.4cm \s^*} \underline{H}( X)(m).}
   \eeq
   We can replace $H(X)(m)$ with, say,  the degree-preserving endomorphisms ${\rm End}^o(H(X))$ 
   of $H(X).$ We thus  obtain  diagram $(*),$ as well as the analogue
   of (\ref{3ef}), for ${\rm End}^o(H(X)),$ and these are the two ingredients needed for the following:
   
      \begin{defi}\la{defah}{\rm
   We say that   $\underline{\zeta} =(\zeta_{dR}, \{\zeta_{et_\ell}\}) \in \underline{H}^{2p}(X)(p)$ is {\em absolute Hodge} if:

\n
1)
   $\s^* \underline{\zeta} \in \im\,\delta_\rat$ (rationality), for every $\s;$

\n
2)
   $\zeta_{dR}
   \in  ({\mathscr F}^0 \cap {\mathscr W}_0 )\, H_{dR}^{2p}(X/K) (p).$

\n
 Similarly, for   $\underline{\zeta} \in \underline{\rm End}^o(H(X)).$ 
   }
   \end{defi}
   
   \begin{rmk}\la{vbg}{\rm
   Membership of $\zeta_{dR}
   \in  ({\mathscr F}^0 \cap {\mathscr W}_0) \, {\rm End}^o(H_{dR}(X/K))$
   means that $\zeta_{dR}$  preserves the Hodge and weight filtrations  on $H_{dR}(X/K).$
   In the Betti context,  membership in  $({\mathscr F}^0 \cap {\mathscr W}_0) \, {\rm End}^o(H_{B}(X,\rat))$ is  equivalent to the  map being of rational  mixed Hodge structures.  
  }\end{rmk}

\subsubsection{Andr\'e motivated endomorphisms}\la{caanmo}
  Let $K$ be an algebraically closed field of characteristic zero and let $X$ be a nonsingular projective $K$-variety.
We freely use Y. Andr\'e's   motivated cycles
$A^{\bullet}_{\rm mot} (X)$ and category  of motives $\m{M}_K$; see \ci{andre}, \S4.  

The $\rat$-vector space 
of motivated cycles of $X$ embeds $\rat$-linearly
into the chosen classical cohomological theory, e.g.: $\xymatrix{cl_H: A^{\bullet}_{\rm mot} (X) 
\ar[r] & H^{2\bullet}_{dR}(X/K)(\bullet)}$.
The category of Andr\'e motives  $\m{M}_K$ --defined the same way as Grothendieck's, but
 using motivated correspondences in lieu of algebraic ones-- is   Tannakian, graded, abelian semi-simple and polarizable. Let $h(X)\in \m{M}_K$ denote the  object (Andr\'e motive)
  associated with a smooth projective $K$-variety $X.$ Since the K\"unneth components of the diagonal are motivated, we have  $h(X)=\oplus_d h^d(X)$. The de Rham  realization
  functors realizes  $h^d(X)$ as 
  $H^d_{dR} (X/K)$,  together with its weight and  Hodge filtrations.
  We have ${\rm End}_{\m{M}_K} (h(X)) = A^{\dim{X}}_{\rm mot} (X \times X)
 \subseteq_{cl_{\underline H}}  H_{dR}^{2\dim{X}} (X\times X/K) (\dim{X}).$
 Note that the last term lives in ${\rm End}^o (H_{dR}(X/K))$.
 
 \begin{defi}\la{motiendo}
 {\rm 
 An endomorphism in ${\rm End}^o (H_{dR}(X/K))$ is {\em Andr\'e motivated}
 if it lies in the image of  $A^{\dim{X}}_{\rm mot} (X \times X)$.
 }
 \end{defi}
Via $cl_H$, a motivated codimension $p$ cycle gives rise to a cycle in
$H^{2p}_{dR}(X/K)(p).$ If $K$ is embeddable in $\comp$, then 
such a cycle gives rise to an   absolute Hodge class (\ci{andre}, Proposition 2.5.1).
Conversely, absolute Hodge cycles are {\em expected} to give rise to  motivated classes.

\subsubsection{Tate classes}\la{tacl}
Let $K_0$ be a field of arbitrary characteristic,  with a fixed  algebraic closure $K_0 \subseteq K$.
Let $X_0$ be a $K_0$-scheme and let $X$ be the resulting $K$-scheme.
For every prime number
$\ell \neq {\rm char} \,K,$ the $\rat_{\ell}$-adic
cohomology (with and without compact supports) of $X$ carries the continuous action of the profinite group
${\rm Gal} (K/K_0)$.

We say that an element  in some $H_{et}(X,\rat_\ell)(m)$, or in ${\rm End}^o (H_{et}(X,\rat_\ell)$,
is {\em of Tate type} if it is ${\rm Gal} (K/K_0)$-invariant.  Similarly, for compact supports.

\subsubsection{Absolutely Hodge and Tate in the sense of Ogus}\la{ahto}
In this section, we  adapt A. Ogus' definitions \ci{ogus} of  {\em absolutely Hodge} and {\em absolutely Tate} de Rham classes to the case of maps. 
Since we work with crystalline cohomology for varieties over a perfect field
of positive characteristic (see the survey
\ci{illusie})  and this theory   is not well-behaved for either singular, or non
proper  varieties, we place ourselves in the safer smooth and proper niche.
Alternatively, one may work with rigid cohomology.

Let  $K$ be a field of characteristic zero,  let $R \subseteq K$ be a smooth $\zed$-algebra,
set
$\m{R}:= {\rm Spec}\,
R$ and let  $\xymatrix{\m{X} \ar[r]& \m{R}}$ be a proper and smooth map. 
It is standard that a proper smooth $X/K$-scheme can be descended to  such
an $\m{X}/R$ and that the closed points of $\m{R}$ have finite, hence perfect, residue fields.
Let $c\in \m{R} (\comp)$ be a $\comp$-point of $\m{R}.$ Denote by $\m{X}_c$ the $\comp$-variety obtained
by base change. We have the natural maps:
\beq\la{1oo}
\xymatrix{
c^*:H_{dR} ({\cal X}/R) \ar[r] & H_{dR}(\m{X}_c/\comp) \ar[r]^\cong & H_B ({\cal X}_c, \comp).}
\eeq
Let $k$ be a perfect field of $\rm{char} \, k = p >0.$
Let $W=W(k) $ be the  ring of Witt vectors of $k.$
Let $\s \in {\cal R}(W)$ be a $W$-point of $\m{R}$ and let 
$\overline{\s} \in {\cal R}(k) $  be the resulting $k$-point, i.e. the  compositum: 
$\xymatrix{R \ar[r] & W \ar[r] &  k}$.
We get the Cartesian diagram:
\beq\la{ghjk}
\xymatrix{
{\cal X}_{\overline{\s}} \ar[r] \ar[d]  & {\cal X}_\s \ar[d]\ar[r] & {\cal X} \ar[d]  \\
 {\rm Spec}\, k  \ar[r] & {\rm Spec}\, W \ar[r] & {\rm Spec} \, R,
}
\eeq
with resulting maps (de Rham/crystalline comparison: $\bb{K}:=$ the quotient field of $W$):
\beq\la{2oo}
\xymatrix{ \overline{\s}^*:
H_{dR} ({\cal X}/R) \ar[r] &
H_{dR} ({\cal X}_\s/ W)  \otimes_W \bb{K} \ar[r]^{\cong_{\bb K}} & H_{\rm cris} ({\cal X}_{\overline{\s}}/W) \otimes_W \bb{K} .
}
\eeq
The {\em crystalline Frobenius} is the automorphism of the rhs  of (\ref{2oo})
induced by the action of 
absolute Frobenius on  $\m{X}_{\overline{\s}}.$ 

\begin{defi}\la{abhott}
{\rm
Let $X$ be smooth and proper over a field $K$ of characteristic zero. Let $\xi \in {\rm End}^o (H_{dR}(X/K)).$ We say that $\xi$ is
{\em absolutely Hodge}, if there are  $\m{X}/R$ as above inducing $X/K,$ and $\xi_R \in
{\rm End}^o (H_{dR}(\m{X}/R))$ inducing $\xi,$  such that,
for every $c\in \m{R}(\comp),$
$c^* \xi_R$ is a Hodge class. We say that $\xi$ is
{\em absolutely Tate}, if there  are 
 $\m{X}/R$ and $\xi_R$  as above, such that,
 for every  $W$ and every $\s\in \m{R}(W),$
$\overline{\s}^* \xi_R$ is invariant under crystalline Frobenius.}
\end{defi}

\subsection{Three geometric constructions and a linear algebra splitting}\la{tgc}

\subsubsection{The filtrations induced by  a flag}\la{fibf}
Let $(K, f)$ be as in (\ref{data}).
Our aim is to define the filtrations (\ref{3ev})   associated with the diagram (\ref{345}).
  Form a Cartesian diagram of $K$-varieties:
   \beq\la{345}
 \xymatrix{
 X'_\bullet \ar[r]  \ar[d]  \ar@/^1pc/[rr]^{r_\bullet}  
 &
 X' \ar[r] \ar[d] & 
 X \ar[d]^f \\
  Y'_\bullet \ar[r]  \ar[d] 
 &
 Y' \ar[r]^q \ar[d]^\iota   & 
 Y \\
 {\Bbb A}_\bullet \ar[r]^{i_\bullet} &
 {\Bbb A},
 }
 \eeq
 where:

 \n 1)
 $q$ is any map subject to the following properties:
 $Y'$ is affine, $q$ is a Zariski locally trivial ${\Bbb A}^d$-fibration for some $d\geq 0;$ the existence of $q,$ which is never unique,
 is ensured by the quasi-projectivity of $Y$ (``Jouanolou's trick" \ci{decpf2}); 
 
\n 2)  $\iota$ is any closed embedding into  an affine space $\Bbb A$
 of some dimension $N;$

 \n
 3)
  $i_\bullet$ is the embedding of a complete  flag of affine linear sections
 of $\bb{A}:$
 
 \n 
 ${\Bbb A}_\bullet = \{ \emptyset =: {\Bbb A}_{-(N+1)}
 \subseteq \ldots \subseteq {\Bbb A}_0:= {\Bbb A}\},$ where 
 ${\Bbb A}_{-c}$ is  a codimension $c$ linear subspace
 of $\Bbb A.$

 \smallskip
The left-most-side of this diagram corresponds
  to a choice of a  point $\beta \in {\cal B} (K),$ where ${\cal B}$ is the variety
  of   complete affine linear flags in ${\Bbb A}.$
  For any such $\beta \in  \m{B}(K)$,
  define  {\em increasing}  filtrations $F^\beta$  on  $H(X)$ and on $H_!(X)$ by setting:
  ($r_{!,k}$ shifts cohomological degrees)
  \beq\la{3ev}
  \xymatrix{
  F^{\beta}_k  H(X) :=  \ke \, \{ H(X)
  \ar[r]^{\hskip -0.9cm r^*_{-k}} & H( X'_{-k} )\}, \;\;
  F^{\beta}_{!,k}  H_!(X) :=  \im \, \{ H_!(X'_k) 
  \ar[r]^{\hskip 2.2cm r_{!, k}} & H_!( X )\}.
  } 
  \eeq
 \begin{fact}\la{pf=fl}{\rm 
 ({\bf Perverse=general flag)}} According to {\rm \ci{decpf2}}, Theorem {\rm 3.3.5}, for {\em general}
 $\beta \in \s\m{B}(\comp),$ the filtrations of type $F^\beta_B$, $F^\beta_{!,B}$ appearing
 at every Betti vertex  of diagram $(*)$ {\rm (\ref{pz1})} for $H(X)$ and $H_!(X)$
 coincide, up to re-numbering, with the 
 perverse filtrations $\ms{P}_{B,f}$.
 \end{fact}

\subsubsection{Teissier and Whitney stratifications of a map}\la{tewhsr}
The goal of this section is to prove Proposition \ref{ratsupp} by using
beautiful results of Verdier \ci{verdier} and Teisssier \ci{teissier}.
For background on stratifications, we refer to \ci{smtt}, p.43.

 Let $(K,f)$ be as in (\ref{data}). If $K =\comp$, then there  are the    notions of:
Whitney stratification of the complex analytic space underlying a $\comp$-scheme;
 Whitney stratification of a map of $\comp$-schemes, which requires the topological local triviality of $f$ over each stratum. Verdier has proved 
 in   \cite{verdier}, 3.3, 4.14-15, that one can produce
 {\em algebraic}  Whitney stratifications,
$\ms{T}_X: X= \coprod_i T'_i$ and $\ms{T}_Y: Y=\coprod_j T_j$,
of the map $f$
so that the strata $T$ are $\comp$-subvarieties  (locally closed, integral and nonsingular).
This is achieved as follows:
Verdier first proves that there are   algebraic  Whitney stratifications  of $Y$ and $X$ (this is the hard part); it is then easy to refine algebraically   both stratifications so that
the following condition is met:
every stratum on $X$ maps smoothly and surjectively onto a stratum
on $Y;$ at this point, $f$ is a stratified submersion for the refined algebraic Whitney
stratifications; the Thom isotopy lemmata \ci{smtt} imply the desired local triviality assertion.

Teissier has introduced a local {\em algebraic} condition
on  a stratification of a $\comp$-schemes that, strikingly, implies that the stratification  is
a Whitney one:  for the condition and the proof of the implication, see  \ci{teissier}, p.379, 
and  Thm. III.2.3.1, p.398.
Teissier's algebraic condition can be defined over any field, where it can be achieved
by a simple noetherian
induction based on the argument of the  proof of 
\cite{teissier}, VI.2.1, p.477.

Given a stratification $\ms{T}_Y: Y=\coprod T_i$, denote by $\m{S}_{\ms{T}_Y}$ the set 
of the closures  $\overline{T_i} \subseteq Y.$

\begin{pr}\label{ratsupp}
Let $(K,f)$ be as in {\rm (\ref{data})}. 
There is a stratification $\ms{T}_Y: Y =\coprod T_j$ 
such that for every embedding $\s$ of $K$ into $\comp$ and for every  $b\in \zed$, the sets of supports {\rm (\S\ref{compbe})}
$\m{S}_B (\s f, b, IC_{B, \s X})\subseteq  \s (\m{S}_{\ms{T}_Y}).$
 In particular the supports of 
every $\s f$ are $K$-rational.
\end{pr}
{\em Sketch of proof.} Start with a stratification  of $Y$ subject to Teissier's
conditions. Do the same for $X$. Refine  both stratifications so that
the following condition is met:
every stratum on $X$ maps smoothly and surjectively onto a stratum
on $Y.$ Call $\ms{T}_Y$ and $\ms{T}_X$ the results, which automatically
satisfy Teissier's conditions.  By passing to $\comp$, the  discussion above shows
that the $\s \ms{T}$ form a 
Whitney stratification of  $\s f$. The local triviality assertion implies that
the direct image complex $\s f_* IC_{B,\s X}$
is constructible with respect to  $\s\ms{T}_Y$, so that
the supports of  $\s f_* IC_{\s X}$  are to be found
in $\s \m{S}_{\ms{T}_Y}$.
\blacksquare

\subsubsection{Contributions of strata: results from \cite{htam, decpf2}}\la{sttrm}
In this section, first we  introduce the  $K$-rational construction (\ref{r111}) (to be used later in the proof of Proposition \ref{hji}), we then use it for $K=\comp$, to recall, 
in the Betti context,  the  linear algebra characterization  (\ref{tsyf})
of the summands corresponding to different supports in the decomposition (\ref{rbb5}).
This characterization follows from 
Claim on p.745 in \cite{htam}  and    Theorem 4.3.2 and Lemma 4.5.3 in \cite{decpf2}.

Let $(K,f)$ be as in (\ref{data}).  Since we are going to use the construction
that follows for $f$, as well for  auxiliary  maps $\xymatrix{\phi : U \ar[r] & V,}$
we use $\phi$ in what follows.
 Fix a stratification $\ms{T}_V$ of $V$ as in Proposition \ref{ratsupp}.
For every stratum $T\subseteq V,$ form  the  commutative diagram:
\beq\la{r111}
 \xymatrix{
 {\cal T} \ar[rr]^{\rho} \ar[rrd]_\gamma   \ar@/^1pc/[rrr]^{r} & & \overline{f^{-1}(\e_T)}_{red} \ar[r] \ar[d] & U \ar[d]^\phi \\
 & & \overline{T} \ar[r] & V,
 }
  \eeq
 where:   $\xymatrix{\overline{T} \ar[r] & V}$ is the  closed  embedding;
$\e_T$ is the generic point of $T;$
  $\xymatrix{\overline{\phi^{-1}(\e_T)}_{red} \ar[r] & U}$ is the closed embedding;
   $\rho$ is a 
  proper, generically finite surjection from a nonsingular  quasi projective $K$-scheme, e.g. first normalize, then resolve
 the irreducible components;  the proper map  $r$  is the evident compositum.
Since $r$ is proper, we have the pull-back maps:
\beq\la{pubama}
\xymatrix{
r^*: H(U) \ar[r] & H(\m{T}), &
r^*_{!}: H_!(U) \ar[r] & H_!(\m{T})}.
\eeq

\smallskip
\n
{\bf Notation.}
In the remainder of this section:  $K=\comp$;  we keep $\ms{T}_V$ and
we assume, in addition that  $U,V$ are quasi projective
varieties with $U$  nonsingular of some dimension $d$; by cohomology we mean Betti cohomology with rational coefficients and we drop the Betti decoration; for $j \in \zed$, define the complementary
degree $j':= 2d-j$ and denote  by $
\xymatrix{\int_{U}: H^j (U) \times H_!^{j'}
(U)  \ar[r] &  \rat}
$
the  Poincar\'e duality pairing; 
finally, set $Gr_{\phi,b}^j:= Gr_{\phi,b}H^j(U)$ and  $Gr_{!,\phi,b}^j:= Gr_{\phi,b}H_!^j(U),$
and similarly for $Gr_{\phi,b, \frak{V}}^j $ and $Gr_{!, \phi,b, \frak{V}}^j.$

\medskip
The supports $\frak{V}$ of $\phi$ are among the  closures  $\overline{T}$
of the strata.  Note that what follows applies to {\em all} the strata $T$,
not just to the ones appearing as supports in a given perversity. This is simply because, in that
case, the corresponding terms are automatically  zero (see \ci{htam}).

 Recalling our conventions on the perverse filtration, the pullback map $r^*: H(U)\lorw H(\cal{T})$ is compatible with the perverse filtrations, provided we shift the one for $\gamma$
 by the fixed amount $d - \dim{\m T}.$ We set, accordingly, $b':= b+ d - \dim{\m T}$. \ Similarly, for the pull-back map $r^*_{!}$ for compact supports.
 
 The map $r^*$ ($r^*_!$, resp.) thus descend to a map $\xymatrix{r^*_b: Gr_{f,b}
 \ar[r] & Gr_{\gamma, b'}}$ ($r^*_{!,b}$, resp.) 
 on the graded pieces, compatibly with the decomposition
 by supports, i.e. $r^*_b= \sum_{\frak{V}} r^*_{b,\frak{V}}$ ($r^*_{!,b}= \sum_{\frak{V}} r^*_{!,b,\frak{V}}$, resp.). Note that the supports for the maps $f$ and $\gamma$ may be quite different, but this does not spoil the picture.
 
Denote by 
$
\xymatrix{r^*_{b,\overline{T}}: Gr_{\phi,b}  \ar[r] & Gr_{\gamma,b', \overline{T}}}$
the compositum of the  map $r^*_b$ 
followed by the projection onto 
the  (possibly trivial) $\overline{T}$-th direct summand of $Gr_{\gamma,b'}.$ Similarly, we have the map   
$\xymatrix{  
{r^*_{!, b,\overline{T}}}:
  Gr_{!, \phi, b}  \ar[r] & Gr_{!,\gamma, b',\overline{T}}
}.$ 

The duality pairing 
descends to a non-degenerate pairing on the graded pieces:
\begin{equation}\label{pdgr}
\xymatrix{\int_{U}^{Gr}:      Gr_{\phi, b}^j       \times Gr_{!,\phi, -b}^{j'} 
  \ar[r] & \rat, & (\alpha,\beta) \ar[r] & \int^{Gr}_U (\alpha,\beta)},
\end{equation}
and the  $Gr^j_{\phi, b, \frak{V}}$ and $Gr^{j'}_{!,\phi, -b,
\frak{V'}}$ are mutually orthogonal when $\frak{V}\neq \frak{V'}.$

We can now state the desired characterization \ci{htam, decpf2}  of the summand corresponding to the {\em dense} stratum $ V^o$: 
\begin{equation}\la{fullsupport}
Gr_{\phi, b, V} = \bigcap_{S \neq V^o  } \ke \, r_{b,\overline{S}}^* \subseteq
Gr_{\phi, b} , \qquad 
Gr_{!,\phi, b, V} = \bigcap_{S \neq V^o  } \ke \, r_{!, b,\overline{S}}^* \subseteq
Gr_{!, \phi, b},
\end{equation}
while, for the {\em non-dense} strata,  the orthogonality property discussed above gives:
\beq \label{tsyf}
Gr_{\phi, b, \overline{T}}=\left(\bigcap_{S \neq T,V^o} \ke \, {r_{b,S}^*}\right) \bigcap \left( Gr_{!,\phi, b, V}  \right) ^{\perp_{\int_U^{Gr}}}.
\eeq
The formula for $Gr_{!,\phi, b, \overline{T}}$ is analogous.

\subsubsection{HL-triples and splitting}\la{hltrspl}
Proposition \ref{00zz} below is one of the keys to the {\em cohomological} approach  
in this paper: it singles out the preferred splitting in cohomology we work with;
see  (\ref{66nn7}) and Proposition \ref{hji}.  It should be compared with 
 \cite{shockwave},  Proposition 2.4, which is stated and proved at the level of {\em triangulated categories
 with $t$-structures}. By taking cohomology, the latter implies the former, in the Betti and \'etale contexts. 
 On the other hand,  Proposition \ref{00zz} becomes useful
 in  contexts  where it is not immediately clear how to use triangulated categories
 (e.g. the cup product operation in the derived category of $D$-modules, in de Rham case), or where
 such formalism is absent (e.g., to our knowledge, crystalline cohomology).

  Let $\cal A$ be an Abelian category endowed with an additive and exact  autoequivalence 
  $\xymatrix{A \ar@{|->}[r] &A(1)},$ whose iterates are denoted by $(m),$ 
  and let ${\cal AF}$ be the associated category of  finitely filtered objects.
  Given $(H,F) \in {\cal AF},$ we denote by $Gr_i H:=Gr^F_i H,$  and by $(Gr_* H,F)$ the associated filtered graded object,
  i.e. $Gr_* H : = \oplus_i Gr_i H$ with the ``direct sum" filtration
  $F_i Gr_* H:= \oplus_{i' \leq i} Gr_iH.$ 
  Let $(H,F,e)$  be a triple with $(H,F)$ in ${\cal AF}$ and $e\!\!: H\xymatrix{ \!\!\ar[r] &H (1)}$
  subject to $\xymatrix{e: F_iH \ar[r] & F_{i+2} H(1)}.$ Assume that the triple
  $(H,F,e)$ is an {\em HL-triple}, i.e. that the  induced
  maps $\xymatrix{e^i: Gr_{-i} H \ar[r] & Gr_i H (i)}$ are isomorphisms (hard Lefschetz property).
  We have the usual primitive Lefschetz decomposition $Gr_*H = \oplus_{0\leq j \leq i}
  P^e_{-i}(-j),$ where $P^e_{-i} := \ke \,e^i  \subseteq Gr_{-i}H.$
  
  \begin{pr}\la{00zz} {\rm (\ci{decasplittings}, Lemma 2.4.1)}
  Let $(H,F,e)$ and $(Gr_*H,F)$ be as above. 
  
  \n
  There is a unique map $\xymatrix{f_i : P^e_{-i} \ar[r] & F_{-i} H}$
  with the following two properties: 
  
  \n
  {\rm 1)} $f_i$ induces the natural inclusion $\xymatrix{P^e_{-i} \ar[r] & Gr_{-i}H};$
  
  \n
  {\rm 2)}
  $\xymatrix{e^{s} \circ f_i : P^e_{-i} \ar[r] & Gr_{s} H (s)}$ is zero $\forall s >i.$
  
  \n
  The map $\xymatrix{\phi_e:=\sum_{0 \leq j \leq i} e^j\circ f_i : (Gr_*H,F) \ar[r] & (H,F)}$ is an isomorphism in
  ${\cal AF}.$
\end{pr}

\subsection{The decomposition theorem and Betti cohomology}\la{dtbet}
Let $(f, g, \eta)$ be as in (\ref{data}), except that we work over $\comp$ and with
rational  Betti cohomology.
Recall the notation in \S\ref{compbe}, especially the one for the graded
pieces of the perverse filtration, and the one in \S\ref{hltrspl} for HL-triples.
  The analogous statements and constructions  for compact supports are  left implicit; e.g.
 write $Gr_{!B,f}$ in (\ref{rbb5}) and $\pi_{!,B} (f,\eta, b)$ in (\ref{poiz}).

\subsubsection{Decompositions associated with a proper map}\la{how11}
If we set $Q:= IC_X$, then we have the finite
direct sum decomposition by supports:
\[
I\!H_B(X) \cong \bigoplus_{b} \bigoplus_{\frak{Y} \in \m{S} (f,b, Q) } Gr_{B, f,b,\frak{Y}}.
\]
The  perverse Leray filtration  ${\mathscr P}_{B,f}$  induced by $f$ on
$I\!H_B(X)$ is given by mixed Hodge substructures
for the natural mixed Hodge structure on $I\!H_B(X)$ constructed in \ci{decpf2},
i.e. it is compatible with $\mathscr{W}$ and, after passing to $\comp$-coefficients,
with $\mathscr{F}.$

 The relative Hard Lefschetz theorem implies that we have direct sum
  isomorphisms
 $\eta_B^b: \oplus_{\frak Y}Gr_{B, f,-b, \frak{Y}} \cong  \oplus_{\frak Y} Gr_{B, f,b, \frak{Y}}(b)$
 of mixed Hodge structures. In particular, we deduce that
 the triple $(I\!H_B( X), {\mathscr P}_{B,f}, \eta_B)$ is an HL-triple (\S\ref{hltrspl})
 for the category of mixed Hodge structures.
We obtain the following  four decompositions of mixed Hodge structures:
\beq\la{rbb5}
Gr_{B,f} := \bigoplus_b Gr_{B, f,b} = \bigoplus_{b,\frak Y} Gr_{B, f, b,{\frak Y}}= 
 \bigoplus_{0\leq j\leq i} P^\eta_{B,-i} (-j)=
\bigoplus_b \bigoplus_{{\frak Y},  -i+2j =b} P^\eta_{B, -i,{\frak Y}} (-j),
\eeq
where $P^\eta_{B, -i, \frak{Y}}(-j): = Gr_{B, f, b,\frak{Y}} \cap P^\eta_{B, -i}(-j).$
We apply Proposition \ref{00zz}  and obtain
the  splitting of mixed Hodge structures of the perverse filtration $\mathscr{P}_{B,f}$: 
\beq\la{io44}
\xymatrix{ \phi_{\eta_B} \, :\; 
Gr_{B,f}
\ar[r]^{\hskip 0.0cm \cong} & I\!H_B(X),
}
\eeq
 no matter which
of the four decompositions 
 (\ref{rbb5})  of $Gr_{B,f}$ we plug
in the l.h.s.

\subsubsection{The projectors   associated with  a projective map}\la{howo}
We have the  four direct sum decompositions stemming from (\ref{rbb5}) and (\ref{io44}):
\beq\la{poi}
I\!H_B(X) \, = \, \bigoplus_b \phi_{\eta_B} \left( Gr_{B,f,b}\right) \, = \, 
\bigoplus_{b,\frak{Y}} \phi_{\eta_B} \left( Gr_{B,f,b,\frak{Y}} \right)\, = \, \ldots.
\eeq
Any direct summand of any direct sum  decomposition above  gives rise to the
degree-preserving  endomorphism  
of $I\!H_B(X)$ obtained by projecting onto that summand. We  thus obtain the {\em 
Betti projectors
of the decomposition theorem for $(f,\eta)$}:
\beq\la{poiz}
\pi_B (f,\eta, b) \;, \pi_B (f,\eta, b, \frak{Y}), \; \pi_B(f, \eta, i,j), \; \pi_B (f, \eta, i,j, {\frak Y})
\; \in    \;  \left(\mathscr{W}_0 \cap \mathscr{F}^0\right) {\rm End}^o (I\!H_B(X)),
\eeq
where the inclusions in the appropriate steps of the filtrations $(\mathscr{W,F})$ express that these projectors
are maps of rational mixed Hodge structures.

  \subsubsection{Decompositions for the composition  of two proper maps}\la{rbrbt}
  We drop the Betti decoration. The purpose of this section is to clarify, via Lemma \ref{44mm} below,
   the relation between
   the perverse filtrations arising from the maps $f,g,h.$
   
  Let $Q$  be semisimple
  complex
  of geometric origin (\ci{bbd}) on $W.$ We have
  the direct sum decompositions
   $g_* Q \cong 
  \bigoplus_a\left(  \bigoplus_{ {\frak X} }  Q_{g,a,{\frak X}}   \right)  [-a]$
  and, since $h_* = f_* \circ g_*$:
   \beq\la{rtm}
  h_* Q \cong \bigoplus_c \left( 
  \bigoplus_{a+b=c} 
  \bigoplus_{\frak X, \frak Y} Q_{ g,a,{\frak X}; f,b,{\frak Y}} 
  \right)
 [-c],
 \eeq
 where it is understood that, for every fixed $(a,b),$ we take,  for each support ${\frak X} \in {\cal S}(g,a, Q),$ the  supports ${\frak Y} \in {\cal S} (f,b, Q_{g,a, \frak{X}}).$

We have the perverse filtrations:
 ${\mathscr P}_g, {\mathscr P}_f$ on $H(W,Q)$ and ${\mathscr P}_f$ on  $Gr_{f,a}:=
 Gr^{{\mathscr P}_f} 
 H(X,Q_{f,a}).$ 
 Clearly, ${\mathscr P}_h$ induces a filtration, still denoted by ${\mathscr P}_h,$ on $Gr_{f,a}.$  
 Given ${\frak X} \in {\cal S} (g,a, Q),$
 consider the natural quotient map $\xymatrix{q:{\mathscr P}_{g,a} \ar[r] & Gr_{g,a}= \oplus_{\frak X}
 Gr_{g,a, {\frak X}} }$ and define ${\mathscr P}_{g,a, {\frak X}} : =q^{-1} Gr_{g,a,{\frak X}}.$
 Clearly, $Gr_{g,a, {\frak X}} = {\mathscr P}_{g,a, {\frak X}}/{\mathscr P}_{g,a-1}.$
 Define ${\mathscr P}_{h,c, {\frak Y}}$ and  ${\mathscr P}_{f,b, {\frak Y}}$ in a similar way.
 
The following lemma  follows directly from the definitions and from (\ref{rtm}).
 
 \begin{lm}\la{44mm} We have the identity:
 \beq\la{p09}
 {\mathscr{P}}_{f,b} Gr_{g,a}  = {\mathscr{P}}_{h,a+b} Gr_{g,a}.
 \eeq
 For every  ${\frak X} \in {\cal S}(g,a,Q)$ and ${\frak Y} \in {\cal S} (f,b, Q_{g,a,{\frak X}}),$
 we have the identity:
 \beq\la{ezrf}
 Gr_{ f,b, {\frak Y} } Gr_{ g,a, {\frak X} } = Gr_{ h, a+b, {\frak Y} } Gr_{ g,a, {\frak X} }
 = \frac{ {\mathscr P}_{h,a+b, {\frak Y}}   \cap { \mathscr{P}_{g,a, {\frak X}}  } }{
 {\mathscr P}_{h,a+b-1}   \cap { \mathscr{P}_{g,a, {\frak X}}  }.
 }
 \eeq
 \end{lm}
 
 \begin{rmk}\la{5566}{\rm 
 The relative hard  Lefschetz theorem applied to $Q_{g,a,\mathscr{X}}$ implies that
 the triple $(Gr_{g,a,{\frak X}}, {\mathscr P}_{f}, \eta)$ is an HL-triple (\S\ref{hltrspl}).
 We deduce the four decompositions of each  $Gr_f\,   Gr_{g,a,\frak{X}}$  as in (\ref{rbb5}),
 as well as the splitting:
 \beq\la{66nn7}
 \xymatrix{ \phi_\eta \, : \; Gr_f\,   Gr_{g,a,\frak{X}} \ar[r] & 
 Gr_{g,a,{\frak X}}
 }
 \eeq  as in (\ref{io44}). When $Q=IC_W,$ everything is  compatible
 with the mixed Hodge structures in sight (\ci{decpf2}).
 We are going to use this set-up in the special case when 
  $g$ is a resolution of the singularities of $X$ integral,  so that $Gr_{g,0,X} = I\!H(X, \rat).$
 }
 \end{rmk}

  \section{The projectors  are absolute Hodge}\la{ciao}

  \subsection{The perverse filtration $\mathscr{P}_f$ is  $K$-rational}\la{bonjour}
  
 \begin{pr}\la{pfg}  Let $(f,\s)$ be as in {\rm (\ref{data})}. 
 There is a unique  filtration ${\mathscr P}_f$ of $H(X),$ which we name the perverse filtration
 associated with $f,$ yielding  diagrams $(*)$ for each subspace 
 ${\mathscr P}_{f,b} H(X)$   and for each inclusion ${\mathscr P}_{f,b} H(X) \subseteq 
 {\mathscr P}_{f,b'}  H(X)$ with $b\leq b'.$ 
 At each Betti vertex,
 the filtration coincides with the Betti perverse filtration $\mathscr{P}_{B,f}.$
 
 \n
 In particular, we can produce diagrams $(*)$ for $\ms{P}_f$ and for 
 $Gr_f:= Gr^{\ms{P}_f}$ on $H(X)$ and on $H_!(X).$
  \end{pr}
  {\em Proof.}
 Given the nature of the arrows of  diagram (\ref{pz1}), unicity  follows from 
   the requirement that, at the Betti vertices,  the filtration coincides with the Betti perverse filtrations.
   
  \n
  As to the existence, we use the construction of \S\ref{fibf}. 
  Since the formation of the maps (\ref{3ev}) is compatible with
  the arrows in (\ref{pz1}), and strictly so for the filtrations $(\mathscr{W,F}),$  the filtrations $F^\beta$
  correspond to each other via the arrows in (\ref{pz1}).
It follows that, for every 
  $\beta \in \m{B}(K)$,  we have  diagrams $(*)$
  for  $H(X)$. Similarly for $F^\beta_!$  and $H_!(X).$
 
  \n
  Recalling that  given a $K$-scheme $Z$, the set $Z(K)$ of its $K$-rational points 
  is dense in every $\s Z (\comp)$, Fact \ref{pf=fl} implies immediately
  that, 
  for $\beta \in {\cal B}(K)$ {\em generic}, the filtrations of type $F^\beta_B$
  and $F^\beta_{!,B}$  appearing in every Betti  
  vertex of (\ref{pz1}) for $H(X)$ and $H_!(X)$ coincide with the  corresponding perverse  filtrations
  $\mathscr{P}_{B,f}.$ 
  \blacksquare

  \subsection{$X$ nonsingular: the supports and the maps  (\ref{rbb5}), 
  (\ref{io44}) are $K$-rational}\la{mi9}
  The goal of this section is to prove
 Proposition \ref{hji}, stating that, when $X$ is {\em nonsingular}, the support and primitive decompositions  (\ref{rbb5}), 
  (\ref{io44}) fit in a diagram $(*).$

\begin{pr}\la{hji}
       Let $(f,\eta, \s)$ be as in {\rm (\ref{data})}. Assume that $X$ is 
        nonsingular.
       There are 
       decompositions {\rm (\ref{rbb5})} and   splittings {\rm (\ref{io44})}
       giving rise to corresponding diagrams $(*)$ in cohomology and in cohomology with compact
       supports.  
\end{pr}
{\em Proof.}   
The proof proceeds by induction on $m:=\dim X.$
We assume that the proposition   holds for every
proper map $\xymatrix{g: Z \ar[r] & Z'}$ of quasi projective varieties
with $Z$ nonsingular and $\dim{Z} < \dim X.$
Note that if $\dim{Z}=0,$ then the inductive hypothesis is
trivial and the conclusion of the proposition is trivially true.

\n
Let $\ms{T}_Y$ be a stratification  of $Y$ as in Proposition  \ref{ratsupp},
so that, for every $\s$,  the supports of the  $\comp$-map $\s f$
are among the closures of the strata of $ \s \mathscr{T}_Y.$

\n
For every  {\em non-dense}  $T\in \mathscr{T}_Y$, form
the $K$-diagram (\ref{r111})   
 for the map $f$.
If necessary, we refine
the stratification so that  it satisfies the conclusion
of Proposition \ref{ratsupp} for  $g$ as well.

\n
Since $\dim{\cal T} < \dim{X}$, we  can apply the inductive hypothesis
to $g$: we have decompositions and splittings (\ref{rbb5}) and (\ref{io44}) for  $Gr_{g,b}$ and $Gr_{!,g,b}$,  giving rise to diagrams $(*)$.
The fact that
$\cal T$ may  fail to be integral does not effect the arguments
in a substantial way. 

\n
Since $r$ is a $K$-map, we have diagram $(*)$ for the  maps
$r^*$ and $r^*_!$ which, 
up to the shift of filtration discussed in \S\ref{sttrm},   are compatible with the respective perverse filtrations:  this is  the case at the Betti vertices, and is thus automatic
at all the other vertices.

\n
By denoting, as in Section \ref{sttrm}, by $\xymatrix{r^*_{b,\overline{T}}: Gr_{f,b} H(X) \ar[r] & Gr_{g,b', \overline{T}} H(\cal{T})}$ 
the map of graded pieces induced by $r^*,$ 
followed by the projection onto the $\overline{T}$-th direct summand of $Gr_{g,b'},$ 
we have diagram $(*)$ for the maps $r^*_{b,\overline{T}}.$ Similarly, for   ${r^*_{!, b,\overline{T}}}.$

\n
Since we have diagram $(*)$ for  the Poincar\'e pairing on $X$,
the usual argument --validity of  a given assertion at the Betti vertex
-- 
coupled with  Proposition  \ref{pfg}, ensures that we have diagram $(*)$
for the graded version  (\ref{pdgr}) of the pairing.

\n
We can therefore {\em define} the summands of the decomposition by supports of $H(X)$ 
at every vertex of diagram $(*)$ by  using the equalities  (\ref{tsyf}) and (\ref{fullsupport})
 in the three cohomology theories.
The resulting decompositions
by supports  give rise to diagrams $(*)$  automatically.

 \n
 Since the primitive decompositions are defined via  the linear algebra
   properties of cupping with $\eta,$ 
  the same is true  for  all four decompositions  in (\ref{rbb5}).

 \n
 Finally, 
 diagram $(*)$ for the splitting  (\ref{io44})  is obtained  formally
 by applying   Proposition \ref{00zz} to each vertex of the diagram.
 \blacksquare

 \subsection{$K$-rational intersection de Rham cohomology}\la{idrc}
 This section can be skipped if in Theorem \ref{maintmz} we assume
  that $X$ is nonsingular.
  
 \begin{pr}\la{rv66}
 Let $X$ be a quasi projective $K$-scheme.
  We have a natural diagram $(*)$ for $I\!H(X)$ and for $I\!H_{!} (X).$
 \end{pr}
 {\em Proof.} We may assume that $X$ is reduced.
 Let $\xymatrix{g: W \ar[r] & X}$ be a ``resolution" of the singularities of
 $X,$ i.e. $W$ is nonsingular and each irreducible component of $W$ is a resolution
 of the singularities of an irreducible component of $X.$ We use the notation of \S\ref{rbrbt}. Note that
 $IC_X$ is the  direct summand $Q_{g, 0,X}$  of $g_* IC_W$ in perversity $0.$

 \n
 Let 
$\mathscr{P}_g$ be   the  perverse filtration  for $g$ on $H(W).$ 
 We apply Proposition \ref{hji} to $g$ and obtain 
 diagram $(*)$ for $Gr_{g,0,X}.$  It only remains to  define the 
 {\em intersection de Rham cohomology groups of $X$} as the $K$-vectorials:
 \beq\la{idrty}
 I\!H_{dR} (X/K) \, := \; Gr_{g, 0,X}  H_{dR}(W/K),
 \eeq
 and similarly for $I\!H_{!,dR}(X/K).$
 \blacksquare

\medskip
It is a routine matter to verify
that these bifiltered groups (they are endowed with the weight and the Hodge filtrations)   are independent of the resolution chosen to define them;
see \ci{htam}, Theorem 2.2.3.a.
 Here is a partial list of useful properties enjoyed by these groups (verifications left to the reader; $X$ integral for simplicity): non-degenerate intersection pairing in intersection cohomology; map
from cohomology to intersection cohomology; 
map $\xymatrix{I\!H_{!,dR} (X/K) \ar[r] & I\! H_{dR} (X/K)};$ 
 cup product with de Rham cohomology $\xymatrix{H^i_{dR} (X/K)
 \otimes I\!H^j_{dR} (X/K) \ar[r] & I\!H^{i+j}_{dR} (X/K)};$ cycle classes; restrictions and Gysin maps
  for normal nonsingular inclusions (\ci{htam}, p. 714]); restrictions to Zariski open subvarieties.

 \subsection{$X$ possibly singular: the maps (\ref{rbb5}) and (\ref{io44}) are $K$-rational}\la{rrm}
  This section can be skipped if in Theorem \ref{maintmz} we assume
  that $X$ is nonsingular.
  
Let things be as  in (\ref{data}).  In the Betti context, 
if we set $Q:= IC_W$, then we have  to the
  decompositions of $Gr_f Gr_{g, a, \frak{X}} H(W)$ and the splitting (\ref{66nn7}) of Remark \ref{5566}.
 
 \begin{pr}\la{rrtt} Let  things be as in {\rm (\ref{data})}  with  $W$ nonsingular. We have diagrams
$(*)$  for the decompositions and splittings
 of Remark {\rm \ref{5566}}

 \end{pr}
{\em Proof.} Since $W$ is nonsingular,  Propositions \ref{pfg} and \ref{hji}
allow us to form diagrams  $(*)$ for the r.h.s. of  the two equalities (\ref{p09}) and
(\ref{ezrf}) of Lemma \ref{44mm}. The conclusion follows.
\blacksquare

\subsection{The projectors of the decomposition theorem are absolute Hodge}\la{99oo}
Let things  be as (\ref{data}), with $g$ a resolution of the singularities of $W.$ 
We apply Proposition \ref{rrtt} in the case of $Gr_{g,0,X},$ which gives us
diagram $(*)$ for $I\!H(X)$ (Proposition \ref{rv66}). Proposition \ref{rrtt}
applies and we obtain diagrams $(*)$ for the decompositions of type
(\ref{rbb5}) of $Gr_{f} Gr_{g,0,X}$ given in Remark \ref{5566},
 and for the splitting $\phi_\eta$ (\ref{66nn7}) in the same remark.
From now on, we write $I\!H(X)$ and $Gr_{g,0,X}$ interchangeably.

We define projectors $\pi (f, \eta, b)$ at each vertex of diagram 
$(*)$ for ${\rm End}^o(I\!H(X))$ by taking the projection
onto the summand $\phi_{\eta} ( Gr_{f,b} Gr_{g,0,X})$ of $I\!H(X)= Gr_{g,0,X}$
associated with $\phi_\eta.$ 
Form diagram $(*)$ for  ${\rm End}^o(I\!H(X) )$
and form the endomorphisms: 
\beq\la{eccoli}
\underline{\pi} (f,\eta, b) \, \in \; \underline{\rm End}^o (I\!H(X)).
\eeq
The collection $\underline{\pi}(f,\eta,b),$ $b \in \zed$ is a  complete system
of orthogonal projectors: 
\[
\underline{\pi}(f,\eta,b) \circ \underline{\pi}(f,\eta,b') =
\delta_{bb'} \,\underline{\pi}(f,\eta,b), \quad \sum_b \underline{\pi}(f,\eta,b) = \underline{{\rm Id}}_{I\!H(X)}.
\]
We can do the same for the remaining
decompositions in (\ref{rbb5}) and obtain projectors $\underline{\pi} (f,\eta, i,j),$ 
$\underline{\pi} ( f, \eta, b ,\frak{Y})$
and $\underline{\pi} (f,\eta, i,j, \frak{Y}).$ These projectors are related
in the following  way: they refine each other in the same way the decompositions
(\ref{rbb5}) refine each other.

We have the analogue $\underline{\pi}_!$ for compact supports.

\begin{tm}\la{maintmz} 
Let $(f, \eta,K)$ be as in {\rm (\ref{data})}. The projectors {\rm (\ref{eccoli})}  
and their compact supports analogue are absolute Hodge.
\end{tm}
{\em Proof.} We have diagram $(*)$ for ${\rm End}^o (I\!H(X)).$ Fix one of the four kind
of projectors  in ${\rm End}^o (I\!H(X))$  and call it $\pi.$ By the very definition of $\pi$ and by Proposition
\ref{rrtt}, we have that the $\pi$'s at each vertex correspond to each other via the
arrows of the diagram. Since the projector is defined  in $\rat$-Betti cohomology,
we have that $\s^* \underline{\pi}$ is rational. Since the projectors are
(strictly) compatible with $(\mathscr{F}, \mathscr{W}),$ we have that $\pi_{dR}
\in (\mathscr{F}^0 \cap \mathscr{W}_0)\, {\rm End}^o (I\!H(X)).$
The compact supports analogue is proved in the same way.
\blacksquare
   
\begin{rmk}
\la{inprod}
{\rm
Assume that $X$ is projective nonsingular. Then Theorem \ref{maintmz},
together  with the K\"unneth formula and Poincar\'e duality, implies that
these projectors    define absolute Hodge classes in  
$\underline{H}^{2 \dim{X}}(X \times X,\rat) (\dim{X}).$  It is not known
 whether these  are
algebraic cycle classes. The Hodge conjecture implies they are.
}
\end{rmk}

  \begin{rmk}\la{mult}
  {\rm
  Let $X_1\stackrel{f_1}\lorw X_2 \stackrel{f_2}\lorw \ldots \stackrel{f_{n-1}}\lorw
  X_n \stackrel{f_n}\lorw$ be proper maps of quasi projective $K$-schemes.
  Set $g_i:= f_i \circ \ldots \circ f_1.$ We leave to the reader the task to formulate and prove 
  by induction on $n$ that there is a
  natural $n$-tuple-version of  the first two  decompositions (\ref{rbb5}),
  involving   multi-supports  of the multi-graded pieces
  $Gr^{g_n}_{b_n} \ldots Gr^{g_1}_{b_1} I\!H(X_1);$ see \S\ref{rbrbt}  for a layout of 
  the case  $n=2.$   There is the further variant, where at each stage
  one takes the primitive decompositions associated with choices of $f_i$-ample line bundles.
The reader can formulate
  and prove the  $n$-tuple variant of 
  Theorem \ref{maintmz} for this  situation.}
  \end{rmk}

\begin{rmk}\la{rty}{\rm
The paper \ci{decasplittings} constructs five distinct distinguished isomorphisms of type (\ref{io44}),
each yielding a collections of projectors as above. Theorem 
\ref{maintmz} remains  valid in  each of these variants and the proofs are similar.}
\end{rmk}

\section{Variants of Theorem \ref{maintmz}}\la{nomen}

\subsection{The projectors are absolute Hodge over  an algebraically 
closed field of characteristic zero}\la{ahcz}
 Theorem \ref{maintmz} is stated for algebraically closed fields with finite transcendence
 degree over the rationals. 
 In \ci{dmos}:
the notion of Absolute Hodge class, or map, is first introduced for  algebraically closed fields
   of characteristic zero with finite transcendence degree over the rationals;
    it is then extended to arbitrary algebraically closed fields
   of characteristic zero by descending to  algebraically closed  subfields of finite
   transcendence degree over the rationals;  finally, it is extended to arbitrary fields   of characteristic zero
  by passing to an algebraic closure  and by considering Galois invariants. 
  It is thus clear that a class, or  map, defined over any algebraically closed field of characteristic
  zero,  is absolute Hodge iff it can be descended to
 an absolute Hodge  class, or map,  over 
  an algebraically closed subfield  of finite transcendence degree  over the rationals.
It is also clear that the  constructions  in  \S\ref{tgc}   are of this nature.

It follows that {\em Theorem {\rm \ref{maintmz}}  holds over any algebraically closed field of characteristic zero and that given a field of characteristic zero, it holds over a suitable finite extension.}

We have defined the groups $I\!H_{dR} (X/K)$   (\ref{idrty})
assuming   $K$  algebraically closed  of finite type over $\rat$. 
What  above shows that one can do so over any field of characteristic zero

 \subsection{The projectors are Andr\'e motivated}\la{am}
  The content of this section is the result of discussions with F. Charles. 
 We place ourselves in the context of \S\ref{caanmo}:
 $\m{M}_K$ is its category of And\'e motives 
 over the algebraically closed field of characteristic zero $K$.  We aim at proving Theorem \ref{tma}. We need a preliminary lemma.

   \begin{lm}\la{adhoc}
  Let $r: X' \lorw X$ be a map of nonsingular $K$-varieties  with $X$ projective.
  There is a natural subobject $k_r \subseteq h(X)$ in $\m{M}_K$ whose 
  de Rham  realization is $\ke \, r^* \subseteq H_{dR}(X/K).$
  \end{lm}
  {\em Proof.}
 The methods of  \ci{ho3} endow de Rham (co)homology with 
  the weight and Hodge  filtrations, and pull-backs and Gysin maps are filtered strict for both. 
  This can be seen by  descending the situation from $K$ to a suitable subfield embeddable into $\comp,$
   by base changing to $\comp,$  and then by using  the bifiltered strict comparison isomorphisms with the Betti theory,
  which enjoys the desired properties. We are thus free to use the usual properties of weights as in
  \ci{ho3}.
  
  \n
  Choose a commutative diagram of  quasi projective $K$-schemes:
  \beq\la{pew}
  \xymatrix{
  & \widetilde{D} \ar[d]^{\iota} \ar[rd]^\nu & \\
  X' \ar[r]^j \ar[rd]_r & \overline{X'} \ar[d]^{\overline r} & \ar[l]_i D \\
  & X,&
  }
  \eeq
  where: $j$ is a smooth projective compactification with simple normal crossing divisor  $D;$ $\nu$ is a resolution of the singularities of $D$ (e.g. the normalization);
  the map $r$ extends to a map $\overline{r}$ (compactify, resolve, put in normal crossing).
  
  \n
  By using the de Rham analogue of  Lefschetz duality, we get, for every $k \geq 0,$  the following commutative diagram  (some decorations omitted):
 \[
  \xymatrix{
  & H^{k-2} (\widetilde{D}) (-1) = H_{2\dim{X'}-k} (\widetilde{D})  (-\dim {X'}) \ar[d]^{\iota_*} \ar[rd]^{\nu_*}&\\
  H^k (X')  & \ar[l]_{j^*} H^k (\overline{X'})  & \ar[l]_{i_*} H_{2\dim{X'} -k }
  (D) (- \dim {X'}) \\
  & \ar[lu]^{r^*} H^k (X) \ar[u]_{\overline{r}^*}, &
  }
  \]
 where:  the
  horizontal line is exact (long exact sequence of relative cohomology); the arrows are filtered strict
  for the Hodge and weight filtrations.
  
  \n
  Clearly, $\ke  \, r^*= {\overline{r}^*}^{-1} (\ke \, j^*)  = {\overline{r}^*}^{-1} 
  (\im \, i_*).$
  
  \smallskip
  \n
  {\bf CLAIM.} We have $\im \, i_* = \im \, \iota_*,$ so that, obviously,
  $\ke \, u^* = {\overline{r}^*}^{-1} 
  (\im \, \iota_*).$
  
  \smallskip
  \n
  {\em Proof of the CLAIM.} By the usual Deligne's mixed Hodge theory,  the map  $\nu_*$ is a surjection
  onto the lowest-weight part ${\ms W}_k H_{2\dim {\cal X} -k} D ( -\dim X').$
  The image of $i_*$ is of pure weight $k$ so that, by strictness, it comes from 
  ${\ms W}_k H_{2\dim {\cal X} -k} D ( - \dim X')$ and the CLAIM follows.
  
 \n
 The category ${\cal M}$ is abelian, so that we can form the compositum:
 \[
 \xymatrix{ v: h(X)  \ar[r]^{\hskip 0.3cm h(\overline{r})} & h(\overline{X'}) \ar[r] & h(\overline{X'})/ \im \, h(\iota).
 }
 \]
 In view of the CLAIM, 
 the kernel of $\ke\, v$  in ${\cal M}$ has Betti realization $\ke \, r^*.$
 
 \n
 Using standard techniques (take a common resolution of two choices), 
 one shows easily that 
$\ke \, v$ is well-defined up to canonical isomorphism independently
of the choices made in the initial commutative diagram. \blacksquare
  
  \begin{tm}\la{tma}
  Let $(f,\eta, K)$ be as in {\rm (\ref{data})} with $K$ algebraically closed of characteristic
  zero and with   $X$   nonsingular and projective.
The projectors of the decomposition theorem  {\rm (\S\ref{99oo})} are motivated cycles in
$A^{\dim{X}}_{\rm mot} (X \times X).$
\end{tm}
{\em Proof.} With Lemma \ref{adhoc} in hand,  the proof of Theorem \ref{tma}  is virtually identical to the one of Theorem
  \ref{maintmz}, which we now review briefly in the context of $\m{M}_K.$
  
  \n
  Lemma \ref{adhoc} and Proposition \ref{pfg} give us motives whose de Rham realization is
  the perverse filtration $\mathscr{P}_f H_{dR}(X).$ We denote these motives
  by $h(\mathscr{P}_f).$ Note the misleading potential of the notation: we are not applying
  $h$ to form this motive, but, rather, we are using Lemma \ref{adhoc} to define it.
  The role of the lemma is precisely to circumvent the fact that in the geometric description of the perverse filtration recalled in Proposition  \ref{pfg},
  one does not immediately ``promote"  diagram (\ref{345}) to the category $\m{M}_K$
  because we are not aware of a theory  of Andr\'e motives for nonsingular open varieties, here the elements of the
  flag $X'_\bullet$ in said diagram. 
  
  \n
  We now produce motives $h (Gr_f),$ decompositions
  of $h(Gr_f)$ as in (\ref{rbb5}), the splitting (\ref{io44}) and, finally,
  motivated projectors $\pi$ in ${\rm End} (h(X)) = A^{\dim{X}}_{\rm mot} (X\times X)$
  which map, via $cl_{dR},$ to our de Rham projectors of the decomposition theorem   in
  $H^{2\dim X}_{dR} (X\times X) (\dim{X}).$
  \blacksquare

\begin{rmk}\la{ght}
{\rm 
The reader should have no problems
  in applying the methods if this paper to define a canonical Andr\'e motive
  attached to the intersection cohomology groups
  of any projective $K$-scheme and prove that Theorem
  \ref{maintmz}, as well as its  variants  outlined in Remarks \ref{mult}, \ref{rty},
   have an intersection
  cohomology  counterpart
  in $\m{M}_K.$
  }
  \end{rmk}

\subsection{The projectors are Tate classes}\la{tate}
Let $(f, \eta)$ be as in (\ref{data}), except that we work over  an arbitrary, i.e. not necessarily algebraically closed field $K.$
The methods of this paper show that
the projectors of the decomposition theorem,
defined  in $\mathscr{W}_0 \, {\rm End}^o
(I\!H_{et} (X\otimes_K \overline{K}, \rat_\ell))$ as in (\ref{poiz}), after passing to an algebraic closure 
$\overline{K}$ of $K,$
are invariant under the action of the group ${\rm Gal}(\overline{K}/K_1),$ where $K_1$ is a suitable finite extension of $K.$ In particular, said projectors are Tate 
for the situation over $K_1$ in the sense of \S\ref{tacl}. 

According to the Tate conjecture,  if $K$ is finitely
generated over it prime field,  and $X$ is geometrically irreducible  smooth and projective,
then these projectors  are expected to be algebraic.

\subsection{The projectors are absolutely Hodge and Tate in the sense of Ogus}\la{rtyog}

\begin{tm}\la{crfr}
Let $(f,\eta)$ be as in {\rm (\ref{data})}, except that we assume that $X,Y$ are projective,
$X$ is smooth and $K$ is a field of characteristic zero. 
After passing to a suitable finite extension of $K,$ the projectors of the decomposition theorem 
are defined and  absolutely Hodge and Tate in the sense of Ogus {\rm (\S\ref{ahto})}.
\end{tm}

\n
{\em Sketch of proof.} The discussion in \S\ref{ahcz} implies that the projectors
are defined after passing to a suitable finite extension of $K.$
We replace $K$ with such an extension.  
The same discussion in \S\ref{ahcz} shows that the projectors are in fact defined over a suitable smooth $\zed$-algebra $R \subseteq K$ over which the whole situation
descends, with $\m{X}/R$ smooth and projective.

\n
Ogus' density argument \cite{ogus}, Remark 4.5 ensures
that in order to prove the ``absolutely Hodge" part of the statement it is enough
to verify it at a $\comp$-point supported at the generic point of $\m{R}.$ This, in turn,
follows immediately from our Theorem \ref{maintmz}.

 \n
 We turn to the ``absolutely
Tate" part of the statement.

\n
As it is recalled in \S\ref{ahcz}, 
the projectors in de Rham cohomology are constructed via diagrams (\ref{345}),
(\ref{r111}) and the splitting $\phi_\eta$ (\ref{66nn7}). In view of the functoriality properties
of crystalline cohomology in the context of smooth and proper varieties,
the use of the second and third ingredient
can be carried out in crystalline cohomology, compatibly with the de Rham/crystalline comparison isomorphism (\ref{2oo}).  On the other hand, the first one involves  non-proper varieties and is thus problematic
on the crystalline side. This issue is easily by-passed via the CLAIM in the proof
of Lemma \ref{adhoc}: one enlarges $R\subseteq K$ by adding finitely many elements so that
diagram (\ref{pew}) descends to it. If necessary, one inverts an element to ensure
smoothness of $\m{X}/R.$ We are now free to form 
the filtration in $H_{\rm cris}(X_{\overline{\s}}/W)$ via the use of diagram
(\ref{345}) base-changed to the perfect field $k,$ compatibly with
the comparison isomorphism (\ref{2oo}). It follows that the formation of every object leading
to the definition of the projectors in crystalline cohomology is now crystalline Frobenius-invariant.
\blacksquare

  \end{document}